\documentclass[leqno,12pt]{amsart}

\usepackage{amsmath,amstext,amssymb,amsopn,amsthm,mathrsfs}
\usepackage{bbm}
\usepackage{verbatim}
\usepackage{enumitem}
\usepackage{color}

\newcommand*\RR{\mathbb{R}}
\newcommand*\NN{\mathbb{N}}
\newcommand*\ind{\mathbbm{1}}
\newcommand*\ven{\vert n\vert}
\newcommand*\al{\alpha}
\newcommand*\fun{\varphi_n}

\newcommand*\Lfun{\mathcal{L}_n^\al}

\newcommand*\Lfunk{\mathcal{L}^{\alpha}_{k}}
\newcommand*\lfun{\ell_n^\al}
\newcommand*\lfuni{\ell^{\alpha_i}_{n_i}}
\newcommand*\lfunk{\ell^{\alpha}_{k}}

\newcommand*\Rop{R_r^{\alpha}}

\newcommand*\Ropi{R_r^{\alpha_i}}

\newcommand*\ve{\varepsilon}

\title[Sharp Hardy's type inequality for Laguerre expansions]
{Sharp Hardy's type inequality for Laguerre expansions}

\author[P{.} Plewa]{Pawe\l{} Plewa}
\address{Pawe\l{} Plewa \newline
			Faculty of Pure and Applied Mathematics, 
      Wroc\l{}aw University of Science and Technology       \newline
      Wyb{.} Wyspia\'nskiego 27,
      50--370 Wroc\l{}aw, Poland      
      }
\email{pawel.plewa@pwr.edu.pl}

\allowdisplaybreaks

\theoremstyle{plain}
\newtheorem{thm}{Theorem}[section]

\newtheorem{lm}[thm]{Lemma}
\newtheorem{prop}[thm]{Proposition}

\theoremstyle{definition}

\theoremstyle{remark}
\newtheorem*{rem*}{Remark}

\theoremstyle{plain}

\begin{document}
\begin{abstract}
A method of proving Hardy's type inequality for orthogonal expansions is presented in a rather general setting. Then sharp multi-dimensional Hardy's inequality associated with the Laguerre functions of convolution type is proved for type index $\al\in[-1/2,\infty)^d$. The case of the standard Laguerre functions is also investigated. Moreover, the sharp analogues of Hardy's type inequality involving $L^1$ norms in place of $H^1$ norms are obtained in both settings.
\end{abstract}

\maketitle
\footnotetext{
\emph{2010 Mathematics Subject Classification:} Primary: 42C10; Secondary: 42B30, 33C45.\\
\emph{Key words and phrases:} Hardy's inequality, Hardy's space, Laguerre expansions of convolution type, standard Laguerre expansions. \\
The paper is a part of author's doctoral thesis written under the supervision of Professor Krzysztof Stempak.
}

\section{Introduction}
Kanjin \cite{Kanjin1} initiated investigation of Hardy's inequalities associated with certain orthogonal expansions. Namely, consider a measure space $(X,\mu)$, where $X$ is one of the domains: $\RR^d$, $(0,\infty)^d$, $d\geq 1$, or $(0,\pi)$, $d=1$, and $\mu$ is the corresponding Lebesgue measure. For a suitable orthonormal basis $\{\fun\}_{n\in\NN^d}$ in $L^2(X,\mu)$, the following inequality was studied
\begin{equation}\label{Hardy_general_introd}
\sum_{n\in\mathbb{N}^d}\frac{\vert \langle f,\fun\rangle\vert}{ (n_1+\ldots+n_d+1)^{E}}\lesssim \Vert f\Vert_{H^1(X,\,\mu)},\qquad f\in H^1(X,\mu),
\end{equation}
where $n=(n_1,\ldots,n_d)$, the symbol $\langle \cdot,\cdot\rangle$ denotes the inner product in $L^2(X,\mu)$, and $H^1(X,\mu)$ is an appropriate Hardy space. The main difficulty lies not only in establishing \eqref{Hardy_general_introd}, but also in finding the smallest admissible exponent $E$, for which such inequality holds. Kanjin's research was inspired by the well known Hardy inequality, which states that (see \cite{HardyLittlewood})
\begin{equation*}
\sum_{k\in\mathbb{Z}}\frac{\vert \hat{f}(k)\vert}{\vert k\vert+1}\lesssim \Vert f\Vert_{\mathrm{ Re}\, H^1}, 
\end{equation*}
where $\hat{f}(k)$ is the $k$-th Fourier coefficient of $f$ and $\mathrm{ Re}\, H^1$ stands for the real Hardy space composed of the boundary values of the real parts of functions in the Hardy space $H^1(\mathbb{D})$, where $\mathbb{D}$ is the unit disk in the plane.

Kanjin \cite{Kanjin1} studied two orthogonal settings: the Hermite and the standard Laguerre functions, both in the one-dimensional case. The obtained admissible exponents were $29/36$ and $1$, respectively. Later Radha \cite{Radha} proved the multi-dimensional version of Hardy's inequality for the Hermite and the special Hermite expansions with appropriate admissible exponents depending on the dimension $d$. The same settings were investigated by Radha and Thangavelu \cite{RadhaThangavelu}. Their result was complemented by Balasubramanian and Radha \cite{BalasRadha}. The one-dimensional version of \eqref{Hardy_general_introd} for the Hermite functions has recently been studied by Z. Li, Y. Yu, and Y. Shi \cite{LiYuShi}. The admissible exponent for the Hermite expansions was finally established as $3d/4$. However, it is not known if the inequality is sharp, that is, if the obtained exponent is the smallest possible.

Hardy's inequality was also considered in the context of other orthogonal expansions. Kanjin and Sato \cite{KanjinSato} investigated the Jacobi setting. Moreover, the case of the multi-dimensional Laguerre expansions of Hermite type was studied by the author in \cite{Plewa}, where the obtained admissible exponent was also equal to $3d/4$. Many authors examined the analogues of \eqref{Hardy_general_introd} replacing the Hardy space $H^1$ by $H^p$ for $p\in(0,1)$ (see \cite{BalasRadha,RadhaThangavelu,Satake}).

Several types of Laguerre function expansions appear in the literature (see e.g. \cite{NowakStempak2, Thangavelu}). In this paper we shall consider two of them: the expansions with respect to the standard Laguerre functions $\Lfun$, and to the Laguerre functions of convolution type $\lfun$. We do not explore the Laguerre polynomial setting because the underlying measure is non-doubling. The studied cases are indeed different in the sense that Hardy's inequality proved in one system does not imply the analogous results for the other expansions. 

The main aim of this article is to investigate Hardy's inequality in the context of the multi-dimensional Laguerre expansions of convolution type, which, up to our knowledge, was never considered before. One of the main novelties is that the measure associated with the underlying space in not Lebesgue measure. However, it is doubling, hence the considered space is a space of homogeneous type (see \cite{CoifmanWeiss}), and therefore the Hardy space is properly defined. 

The next new important feature of the obtained inequality is its sharpness. The explicit counterexample is constructed to show that the admissible exponent cannot be lowered. It turns out that in both Laguerre-type expansions it is possible to point out very similar counterexamples. Therefore, we devote Section \ref{S4} to the standard Laguerre functions. Since the result for the functions $\Lfun$ is known in the literature only in the one-dimensional setting, and the method described in Section \ref{S2} is easily applicable in this situation, we present the proof of appropriate Hardy's inequality for an arbitrary dimension.

%The range of the Laguerre type multi-index $\al$ is the set $[-1/2,\infty)^d$ for the Laguerre functions of convolution type and the set $[0,\infty)^d$ for the standard Laguerre functions. In the second case we needed to fill an emergent gap in the admissible range of $\al$'s (see Lemma \ref{Lfun_alfa_lem}). 

In \cite{Kanjin2} Kanjin proposed investigation of the $L^1$-analogues of \eqref{Hardy_general_introd} in which one replaces $H^1$ norm by $L^1$ norm. He studied such an inequality in the contexts of the Hermite function expansions and the standard Laguerre function expansions, both for $d=1$. This issue was also explored in \cite{Plewa} for the multi-dimensional Laguerre expansions of Hermite type. In this paper we also consider this type of inequality in the Laguerre setting of convolution type and in the standard Laguerre setting. Therefore, Sections \ref{S3} and \ref{S4} are both divided onto two subsections: one deals with a version of \eqref{Hardy_general_introd} and the second with its $L^1$-analogue. The method of proving Hardy's inequality is described in Section \ref{S2}, and relies on an idea proposed in the one-dimensional version in \cite{LiYuShi} and developed by the author in \cite{Plewa}. Roughly, it consists in estimating the derivatives of kernels of certain family of integral operators associated with the considered orthogonal basis. On the other hand, the main tools in establishing $L^1$-type inequalities are the pointwise asymptotic estimates for the functions composing the investigated expansions. 

We shall frequently use, without any further mention, the two basic estimates: for $t,\,T>0$ we have $\sup_{x>0} x^t e^{-Tx}<\infty$, and $(n_1+\ldots+n_d+1)^d\geq (n_1+1)\cdot\ldots\cdot(n_d+1)$, where $n_i\in\mathbb{N}$, $i=1,\ldots,d$. 

\subsection*{Notation}
Throughout this paper the symbol $d\geq 1$ denotes the dimension. We shall use the notation $\RR^d_+=(0,\infty)^d$ and $\NN_+=\NN\setminus\{0\}=\{1,2,\ldots\}$. We write $u,v$ for real one-dimensional variables, $k$ or $j$ for non-negative integers, and $x=(x_1,\ldots,x_d),\ y=(y_1,\ldots,y_d)$ for real multi-dimensional variables. The Euclidean norm is denoted by $|x|$ and $|y|$. Similarly, $n=(n_1,\ldots,n_d)\in\NN^d$  shall stand for a multi-index and $\ven=n_1+\ldots+n_d$ for its length. For the constant multi-indices we will use the bold font, e.g. $\bold{0}=(0,\ldots,0)$. The Laguerre type multi-index $\al=(\al_1,\ldots,\al_d)\in(-1,\infty)^d$ will be denoted by the same symbol $\al$ for $d=1$ and $d\geq 1$, but we hope that it will be always clear from the context whether $\al$ refers to $d=1$ or $d\geq 1$. Again, $|\al|=\al_1+\ldots+\al_d$, stands for the length of the multi-index $\al$. Note that $|\al|$ may be negative. We will use the usual convention writing $x^\al=\prod_{i=1}^d x_i^{\al_i}$. For functions $f,\,g\in L^2(X,\mu)$, where $(X,\mu)$ is a measure space, we denote the standard inner product by $\langle f,g\rangle$. The space $(X,\mu)$ may differ between sections, so the inner products differ as well. Nevertheless, we will use the same symbol throughout the whole paper. If the measure $\mu$ is Lebesgue measure, then we will simply write $L^2(X)$ (and analogously for other function spaces).

We shall use the symbol $\lesssim$ denoting an inequality with a constant that may depend on the parameters that appear before the inequality, but does not depend on the ones quantified afterwards. Also, the symbol $\simeq$ means that $\lesssim$ and $\gtrsim$ hold simultaneously.

\subsection*{Acknowledgement}
The author would like to thank Professor Krzysztof Stempak for his valuable remarks leading to an improvement of the presentation.

\section{General setting}\label{S2}
In this section we present an improved version of the method introduced in \cite{LiYuShi} and developed in \cite{Plewa}. In the next sections we shall apply the result in specific settings.

Let $(X,\mu)$ be a measure metric space such that $X$ is an open convex subset of $\RR^d$, the measure $\mu$ is doubling, and the space is equipped with the Euclidean metric. Note that this implies that $(X,\mu)$ is a space of homogeneous type in the sense of Coifman and Weiss (see \cite[pp.~587-588]{CoifmanWeiss}).

Moreover, let $\{\fun\}_{n\in\NN^d}$ be an orthonormal basis in $L^2(X,\mu)$. We define the family of operators $\{R_r\}_{r\in(0,1)}$ via
\begin{equation}\label{R_def}
R_r f=\sum_{n\in\NN^d} r^{\ven} \langle f,\fun\rangle \fun,\qquad r\in(0,1),\qquad f\in L^2(X,\mu),
\end{equation}
where $\langle\cdot,\cdot\rangle$ denotes the inner product in $L^2(X,\mu)$. Notice that for every $r\in(0,1)$ the operator $R_r$ is a contraction on $L^2(X,\mu)$.

We impose the following assumptions:
\begin{enumerate}[label=(A\arabic*)]
\item \label{A1} $\fun\in L^{\infty}(X,\mu)$, $n\in\NN^d$;
\item \label{A2} there exists $N$ such that
\begin{equation*}
\mu(B(x,\rho))\gtrsim \rho^N,
\end{equation*}
uniformly in $x\in X$ and $\rho\in (0,\mathrm{diam}(X))$, where $B(x,\rho)=\{y\in X\colon \sum_{i=1}^d (x_i-y_i)^2<\rho^2\}$;
\item \label{A3} the operators $R_r$ are integral operators and the associated kernels satisfy for some $\gamma>0$ and a finite set $\Delta$ composed of positive numbers the condition
\begin{equation*}
\Vert R_r(x,\cdot)-R_r(x',\cdot) \Vert_{L^2(X,\,\mu)}\lesssim \sum_{\delta\in\Delta}|x-x'|^\delta(1-r)^{-\frac{\gamma(N+2\delta)}{N+2}},
\end{equation*}
uniformly in $r\in(0,1)$, $x'\in X$, $\rho\in(0,\mathrm{diam}(X))$ and almost every $x\in B(x',\rho)$, where $x'$ and $\rho$ are such that $\mu(B(x',\rho))\leq 1/2$ (we remark that $1/2$ is chosen arbitrarily; in some settings it is convenient to take a smaller positive constant).
\end{enumerate}

The parameter $N$ in \ref{A2} is equal $d$ for Lebesgue measure and, on the whole, is always not smaller than $d$. The smaller the parameter $N$ is, the better for our purposes (however, in general, the smallest $N$ satisfying \ref{A2} may not exist).

Note that that \ref{A3} (with $\Delta=\{1\}$) is implied by much easier condition, namely the operators $R_r$ are integral operators, the kernels $R_r(x,y)$ are differentiable almost everywhere with respect to the first variable, and for $\gamma$ there is
\begin{equation}\label{cond C}\tag{A3'}
\mathrm{ess\,sup}_{x\in X} \big\Vert |\nabla_x R_r(x,\cdot)| \big\Vert_{L^2(X,\,\mu)}\lesssim (1-r)^{-\gamma},\qquad r\in(0,1).
\end{equation}

The $(1,2)$-atoms in the sense of Coifman and Weiss (see \cite[p.~591]{CoifmanWeiss}), which in this paper are called $H^1(X,\mu)$-atoms, are measurable functions $a$ supported in balls $B(x_0,\rho)$, $x_0\in X$, $\rho\in(0,\mathrm{diam}(X))$, such that
$$\int_B a(x)\,d\mu(x)=0,\qquad \Vert a\Vert_{L^2(X,\,\mu)}\leq \mu(B)^{-1/2}.$$
We stress that if $\mu(X)<\infty$, then, additionally, $a\equiv\mu(X)^{-1/2}$ is also considered as an $H^1(X,\mu)$-atom. The (atomic) Hardy space $H^1(X,\mu)$ (compare \cite[pp.~591-592]{CoifmanWeiss}) is composed of functions $f\in L^1(X,\mu)$ admitting the atomic decomposition
\begin{equation}\label{atomic_decomposition}
f=\sum_{j=0}^\infty \lambda_j a_j,
\end{equation}
where $a_j$'s are $H^1(X,\mu)$-atoms, $\sum_{j=0}^\infty|\lambda_j|<\infty$, and the series in \eqref{atomic_decomposition} is convergent in $L^1(X,\mu)$. Moreover, $H^1(X,\mu)$ is a Banach space with the norm
$$\Vert f\Vert_{H^1(X,\,\mu)}=\inf\sum_{j=0}^\infty|\lambda_j|,  $$
where the infimum is taken over all representations as in \eqref{atomic_decomposition}. It is worth mentioning that for every $f\in H^1(X,\mu)$ there is
\begin{equation}\label{H1_L1_norm_ineq}
\Vert f\Vert_{L^1(X,\,\mu)}\leq \Vert f\Vert_{H^1(X,\,\mu)}.
\end{equation}

We emphasize that although we use the $(1,2)$-atoms instead of the usual $(1,\infty)$-atoms, the main results are also valid for the atomic Hardy spaces based on the latter, because the implied Hardy spaces coincide and the associated norms are equivalent (see \cite[p.~592]{CoifmanWeiss}). 

The following lemma holds.

{\lm\label{Lemma_S2} If the assumptions \ref{A1}-\ref{A3} are satisfied, then
$$\int_0^1\Vert R_r a\Vert_{L^2(X,\,\mu)}(1-r)^{(\gamma N)/(N+2)-1}\,dr\lesssim 1, $$
uniformly in $H^1(X,\mu)$-atoms.}
\begin{proof}
If $\mu(X)<\infty$ and $a\equiv\mu(X)^{-1/2}$, then the claim holds trivially. Hence, let us fix an $H^1(X,\mu)$-atom $a$ supported in a ball $B$ such that $\int_B a(x)\,d\mu(x)=0$. Let $x'\in X$ be the center of $B$. Note that since $R_r$ are contractions on $L^2(X,\mu)$ we have for every $0<r<1$
$$\Vert \Rop a\Vert_{L^2(X,\,\mu)}\leq\Vert a\Vert_{L^2(X,\,\mu)}\leq \mu(B)^{-1/2}.$$
This gives the claim in the case $\mu(B)\geq 1/2$. From now on, let us assume that $\mu(B)<1/2$. Minkowski's integral inequality, \ref{A3}, \ref{A2}, and H\"{o}lder's inequality imply
\begin{align*}
\Vert R_r a\Vert_{L^2(X,\,\mu)}&=\Big( \int_{X}\Big\vert\int_{B}\big(R_r(x,y)-R_r(x',y)\big)a(x)d\mu(x)\Big\vert^2  d\mu(y)\Big)^{1/2}\\
&\lesssim\int_{B}|a(x)|\sum_{\delta\in\Delta}|x-x'|^\delta(1-r)^{-\frac{\gamma(N+2\delta)}{N+2}}  d\mu(x)\\
&\lesssim \sum_{\delta\in\Delta}\mu(B)^{\delta/N} (1-r)^{-\frac{\gamma(N+2\delta)}{N+2}}.
\end{align*}
Thus, using the above estimates we obtain
\begin{align*}
\int_0^1 \Vert R_r a\Vert_{L^2(X,\,\mu)}(1-r)^{\gamma N/(N+2)-1}dr&\lesssim \sum_{\delta\in\Delta}\int_0^{1-\mu(B)^{(N+2)/2N\gamma}} \mu(B)^{\delta/N} (1-r)^{-\frac{2\gamma\delta}{N+2}-1}dr\\
&\quad+\int_{1-\mu(B)^{(N+2)/2N\gamma}}^1 \mu(B)^{-1/2}(1-r)^{\gamma N/(N+2)-1}dr,
\end{align*}
and this quantity is bounded by a constant that does not depend on $\mu(B)$.
\end{proof}

{\thm\label{H1_general} Assume that \ref{A1}-\ref{A3} are satisfied. The inequality 
$$\sum\limits_{n\in\NN^d}\frac{|\langle f,\fun\rangle|}{(\ven+1)^E}\lesssim \Vert f\Vert_{H^1(X,\,\mu)}, $$
holds uniformly in $f\in H^1(X,\mu)$, where
\begin{equation}\label{Exponent_formula}
E=\frac{\gamma N}{(N+2)}+\frac{d}{2}.
\end{equation}
} 

The proof of Theorem \ref{H1_general} is almost identical to the proof of \cite[Theorem~4.2]{Plewa} and could be skipped; for the reader's convenience we present it in the Appendix.

We remark that instead of imposing \ref{A2} and \ref{A3}, we could assume that a version of Lemma \ref{Lemma_S2} holds, namely for certain $\zeta>0$
$$\int_0^1\Vert R_r a\Vert_{L^2(X,\,\mu)}(1-r)^{\zeta-1}\,dr\lesssim 1,$$
uniformly in $H^1(X,\mu)$-atoms. Then the implied admissible exponent in Theorem \ref{H1_general} would be $E=\zeta+d/2$.

\section{Laguerre setting of convolution type}\label{S3}
The one-dimensional {\it Laguerre functions of convolution type} of order $\al>-1$ on $\RR_+$ are the functions
$$\lfunk(u)=\Big(\frac{2\Gamma(k+1)}{\Gamma(k+\al+1)}  \Big)^{1/2}L_{k}^{\al}(u^2)e^{-u^2/2},\qquad u>0.$$
In higher dimension the functions $\lfun(x)$ are defined as tensor products of the one-dimensional functions $\lfunk$. Moreover, the Laguerre functions of convolution type form an orthonormal basis in $L^2(\RR_+^d,d\mu_\al)$, where $d\mu_\al(x)=x^{2\al+\bold{1}}\,dx$. The measure $\mu_\al$ is doubling for $\al\in [-1/2,\infty)^d$ (see e.g. \cite[Appendix~1]{AnkerDziubanski}).

We shall verify assumptions \ref{A1}-\ref{A3}. Firstly, we will use the pointwise asymptotic estimates (compare \cite[p.~435]{Muckenhoupt} and \cite[p.~699]{AskeyWainger})
\begin{equation}\label{lfun_asympt}
\vert \lfunk(u)\vert\lesssim\left\{ \begin{array}{ll}
\nu^{\al/2},  & 0<u\leq \nu^{-1/2},\\
u^{-\al-1/2}\nu^{-1/4}, & \nu^{-1/2}<u\leq\sqrt{\nu/2},\\
u^{-\al}(\nu(\nu^{1/3}+\vert u^2- \nu\vert))^{-1/4}, & \sqrt{\nu/2}<u\leq \sqrt{3\nu/2},\\
u^{-\al}\exp(-\gamma u), & \sqrt{3\nu/2}<u<\infty,
\end{array}\right.
\end{equation}
where $\nu=\nu(\al,k)=\max(4k+2\al+2,2)$ and $\gamma>0$ depends only on $\alpha$.

There following formula for the derivatives of functions $\lfunk$,
\begin{equation}\label{lfun_derivative_formula}
\frac{d}{d u}\lfunk(u)=-2\sqrt{k}u\ell_{k-1}^{\al+1}(u)-u\lfunk(u),
\end{equation}
where $\ell_{-1}^{\al}\equiv 0$, is known. Combining \eqref{lfun_asympt} and \eqref{lfun_derivative_formula} we get for $\al>-1$,
\begin{equation}\label{lfun_boundedness}
\Vert \lfunk\Vert_{L^\infty(\RR_+)}\lesssim (1+k)^{|\frac{\alpha}{2}+\frac{1}{6}|-\frac{1}{6}},\qquad \Vert(\lfunk)'\Vert_{L^\infty(\RR_+)}\lesssim (1+k)^{|\frac{\al}{2}+\frac{7}{24}|+\frac{11}{24}}.
\end{equation}
For the first estimate compare \cite[p.~87]{StempakTohoku}.

\subsection{$H^1$ result for Laguerre functions of convolution type}
We consider the family of operators $\{\Rop\}_{r\in(0,1)}$ defined in \eqref{R_def}, associated with the functions $\lfun$. The corresponding kernels are given by
$$\Rop (x,y)=\sum_{n\in\NN^d} r^{\ven} \lfun(x)\lfun(y)=\prod_{i=1}^d \sum_{n_i=0}^\infty r^{n_i} \lfuni(x_i)\lfuni(y_i),\qquad x,y\in\RR^d_+,\ r\in(0,1). $$
An explicit formula for the one-dimensional kernels is known, namely (compare \cite[p.~102]{Szego})
\begin{equation}\label{lfun_Rker_explic}
\Rop(u,v)=\frac{2 }{(1-r)r^{\al/2}(uv)^{\al}}\exp\Big(-\frac{1}{2}\frac{1+r}{1-r}(u^2+v^2)\Big)I_{\al }\Big(\frac{2r^{1/2}}{1-r}uv\Big),\quad u,v>0,
\end{equation}
where $I_{\alpha}$ denotes the modified Bessel function of the first kind, which is smooth and positive on $(0,\infty)$.

The following asymptotic estimates (see \cite[p.~136]{Lebedev}) for the Bessel function $I_s$ are known
\begin{equation}\label{Bessel_estim}
I_{s}(u) \lesssim u^{s}, \quad 0<u<1,\qquad I_{s}(u) \lesssim u^{-1/2}e^{u}, \quad u\geq 1,
\end{equation}
where $s> -1$.

We shall estimate the kernels $\Rop(u,v)$. Note that \eqref{Bessel_estim} gives
\begin{equation}\label{lfun_Rker_estim}
\Rop(u,v)\lesssim\left\{\begin{array}{ll}
(1-r)^{-\al-1}\exp\left(-\frac{1}{2}\frac{1+r}{1-r}(u^2+v^2) \right), & v\leq \frac{1-r}{2\sqrt{r}u},\\
(1-r)^{-1/2}r^{-\al/2-1/4}(uv)^{-\al-1/2}\exp\left(-\frac{1}{2}\frac{1+r}{1-r}(v-u)^2-\frac{uv(1-r)}{(1+\sqrt{r})^2} \right), & v\geq \frac{1-r}{2\sqrt{r}u}.
\end{array}\right.
\end{equation}

We remark that the computation in the proofs of Lemma \ref{lfun_Rker_lemma} and Proposition \ref{lfun_derRker_prop}, as well as in the proofs of Lemma \ref{Lfun_Rker_standard} and Proposition \ref{Lfun_derRker_prop}, are uniform in the indicated ranges of $r$ and $u$. 

{\lm\label{lfun_Rker_lemma} For $\al\geq -1/2$ there is
$$\sup_{u\in\RR_+}\Vert \Rop(u,\boldsymbol{\cdot})\Vert_{L^2(\RR_+,\,d\mu_\al)}\lesssim (1-r)^{-(\al+1)/2},\qquad r\in(0,1). $$}
\begin{proof}
For $0<r\leq 1/2$, we use Parseval's identity and \eqref{lfun_boundedness} obtaining
$$\sup_{u\in\RR_+}\Vert \Rop(u,\boldsymbol{\cdot})\Vert_{L^2(\RR_+,\,d\mu_\al)}\leq  \Big(\sum_{k=0}^\infty 2^{- 2k}\Vert \lfunk\Vert_{L^{\infty}(\RR_+)}^2\Big)^{1/2}\lesssim 1.$$

For $1/2<r<1$, we denote $v_0=(1-r)/(2\sqrt{r}u)$, and estimate the involved integrals over $(0,v_0]$ and $(v_0,\infty)$. Thus, the substitution $\eta=(v\sqrt{1+r})/\sqrt{1-r}$ gives for the first integral
\begin{align*}
\int_0^{v_0} \Rop(u,v)^2\,d\mu_\al(v)&\lesssim  (1-r)^{-2\al-2}\int_0^{v_0}\exp\Big(-\frac{1+r}{1-r}v^2\Big)v^{2\al+1}\,dv\\
&\lesssim (1-r)^{-2\al-3/2}\int_0^{\infty}(1-r)^{\al+1/2}\eta^{2\al+1}e^{-\eta^2}\,d\eta\\
&\lesssim (1-r)^{-\al-1}.
\end{align*}
 Similarly, for the second integral, for $u\leq v_0$ we have
\begin{align*}
\int_{v_0}^{\infty} \Rop(u,v)^2\,d\mu_\al(v)&\lesssim  (1-r)^{-2\al-2}\int_{v_0-u}^{\infty}\exp\Big(-\frac{1+r}{1-r}v^2\Big)(u+v)^{2\al+1}\,dv\\
&\lesssim  (1-r)^{-2\al-3/2}\int_{0}^{\infty}(1-r)^{\al+1/2}(1+\eta^{2\al+1})e^{-\eta^2}\,d\eta\\
&\lesssim (1-r)^{-\al-1},
\end{align*}
and for $u\geq v_0$
\begin{align*}
\int_{v_0}^{\infty} \Rop(u,v)^2\,d\mu_\al(v)&\lesssim  (1-r)^{-1}\int_{v_0-u}^{\infty}u^{-2\al-1}\exp\Big(-\frac{1+r}{1-r}v^2\Big)\,dv\\
&\lesssim  (1-r)^{-1}\int_{-\infty}^{\infty}(1-r)^{-\al-1/2}\exp\Big(-\frac{1+r}{1-r}v^2\Big)\,dv\\
&\lesssim  (1-r)^{-\al-1}\int_{-\infty}^{\infty}e^{-\eta^2}\,d\eta\\
&\lesssim (1-r)^{-\al-1}.
\end{align*}
Combining the above gives the claim.
\end{proof}

A formula for the derivative of the Bessel function $I_\al$ (see \cite[p.~110]{Lebedev}) gives

\begin{align}
&\partial_u\Rop(u,v)=\frac{2ruv^2}{1-r}R^{\alpha+1}_r(u,v)-\frac{(1+r)u}{1-r}\Rop(u,v)\nonumber\\
&=\frac{2(uv)^{-\al}}{(1-r)^2 r^{\al/2}}\exp\Big(-\frac{(1+r)(u^2+v^2)}{2(1-r)}\Big)\Big(2\sqrt{r}v(I_{\al+1}-I_\al)+\big(2\sqrt{r}v-(1+r)u\big)I_\al\Big),\label{lfun_derKer_long}
\end{align}
where all of the above Bessel functions are taken in $\frac{2\sqrt{r}uv}{1-r}$.

Hence, applying \eqref{Bessel_estim} and the first identity from \eqref{lfun_derKer_long} we get the estimates
\begin{align}\label{lfun_derRker_estim}
&\Big|\partial_u\Rop(u,v)\Big|^2\nonumber\\
&\lesssim\left\{\begin{array}{ll}
(1-r)^{-2\al-4}\exp\left(-\frac{1+r}{1-r}(u^2+v^2) \right)(u^2+v^2), & v\leq \frac{1-r}{2\sqrt{r}u},\\
(1-r)^{-3}(\sqrt{r}uv)^{-2\al-1}\exp\left(-\frac{1+r}{1-r}(v-u)^2-\frac{2uv(1-r)}{(1+\sqrt{r})^2} \right)(u^2+v^2), & v\geq \frac{1-r}{2\sqrt{r}u}.
\end{array}\right.
\end{align}

{\lm\label{Bessel_ratio_estim} For $\alpha\geq -1/2$ there is
$$\left\vert I_{\al+1}(u)-I_{\al}(u)\right\vert\lesssim u^{-1} I_{\al+1}(u),\qquad u>0. $$ }
For the proof see \cite[pp.~6-7]{Nasell}.

Applying Lemma \ref{Bessel_ratio_estim} to the second identity from \eqref{lfun_derKer_long} we obtain the estimate
\begin{equation}\label{lfun_derKer_cancellation_estim}
|\partial_u \Rop(u,v)|^2\lesssim v^2|R_r^{\al+1}(u,v)|^2+(2\sqrt{r}v-(1+r)u)^2|\Rop(u,v)|^2(1-r)^{-2} .
\end{equation}

%\begin{equation}\label{Bessel_explicit_formulas}
%I_{-1/2}(u)=\Big(\frac{2}{\pi u}\Big)^{1/2}\cosh u,\qquad I_{1/2}(u)=\Big(\frac{2}{\pi u}\Big)^{1/2}\sinh u,\qquad u>0.
%\end{equation}

%%%Note that Lemma \ref{deri_R_ker} works for $\al>0$, but we want to include the case $\al=-1/2$ as well. Thus, using \eqref{def_R_ker_explic} and \eqref{explicit_formulas} we obtain
%%%$$\Ropal(x,y)=\frac{2}{\sqrt{\pi}}(1-r)^{-1/2}\exp\Big( -\frac{1}{2}\frac{1+r}{1-r}(x^2+y^2)\Big) \cosh\Big(\frac{2\sqrt{r}xy}{1-r}\Big). $$
%%%Hence,
%%%\begin{align}\label{ker_-1/2_deri}
%%%\big(\partial_x\Ropal (x,y)\big)^2&=\frac{4}{\pi}(1-r)^{-3}\exp\Big(-\frac{1+r}{1-r}(x^2+y^2)\Big)\nonumber\\
%%%&\times \bigg( 2\sqrt{r}y\sinh\Big(\frac{2\sqrt{r}xy}{1-r}\Big)-(1+r)x\cosh\Big(\frac{2\sqrt{r}xy}{1-r}\Big)\bigg)^2.
%%%\end{align}
%%
%%%Using basic estimates for $\cosh$ and $\sinh$ and combining \eqref{ker_-1/2_deri} with \eqref{R_ker_estim} and Lemma \ref{deri_R_ker} we obtain for $\al\in\{-1/2\}\cup[1/2,\infty)$
%%%\begin{align}\label{final_estimate_without_cancellations}
%%%&\big(\partial_x\Rop (x,y)\big)^2\nonumber\\
%%%&\lesssim\left\{\begin{array}{ll}
%%%(1-r)^{-2\al-2}(xy)^{2\al+1}(A_\al (x,y)+x^2(1-r)^{-2})\exp\left(-\frac{1+r}{1-r}(x^2+y^2) \right), & y\leq \frac{1-r}{2\sqrt{r}x},\\
%%%(1-r)^{-3}(x^2+y^2+(1-r)^2x^{-2})\exp\left(-\frac{1+r}{1-r}(y-x)^2-\frac{2xy(1-r)}{(1+\sqrt{r})^2} \right), & y\geq \frac{1-r}{2\sqrt{r}x},
%%%\end{array}\right.
%%%\end{align}
%%%where $A_\al(x,y)=x^{-2}$ for $\al\geq 1/2$ and $A_{-1/2}(x,y)=y^2(1-r)^{-2}$.

{\prop\label{lfun_derRker_prop} For $\al\geq -1/2$ we have
$$\sup_{u\in\RR_+}\left\Vert\partial_u\Rop(u,\boldsymbol{\cdot})\right\Vert_{L^2(\RR_+,\,d\mu_\al)}\lesssim (1-r)^{-(\al+2)/2},\qquad r\in(0,1). $$}
\begin{proof}
If $0<r\leq 1/2$, then Parseval's identity and \eqref{lfun_boundedness} yield
\begin{equation*}
\sup_{u>0}\left\Vert\partial_u\Rop(u,\boldsymbol{\cdot})\right\Vert_{L^2(\RR_+,\,d\mu_\al)}\leq\Big(\sum_{k=0}^\infty 2^{-2k}\Vert(\lfunk)'\Vert_{L^{\infty}(\RR_+,\,d\mu_\al)}^2\Big)^{1/2}\lesssim 1.
\end{equation*}

From now on we assume that $1/2<r<1$. We use the notation $v_0=(1-r)/(2\sqrt{r}u)$ again, and split the integration over two intervals: $(0,v_0]$ and $(v_0,\infty)$. In the first case, \eqref{lfun_derRker_estim} and the substitution $v=\eta\sqrt{1-r}$ give
\begin{align*}
\int_0^{v_0}\left(\partial_u\Rop(u,v)\right)^2\,d\mu_\al(v)&\lesssim (1-r)^{-2\al-3}\int_0^{v_0}\exp{\left(-\frac{v^2}{1-r}\right)}v^{2\al+1}\, dv\\
&\lesssim (1-r)^{-\al-2}\int_0^{\infty}\eta^{2\al+1} e^{-\eta^2}\, d\eta\\
&\lesssim (1-r)^{-\al-2}.
\end{align*} 

The case of integration over $(v_0,\infty)$ is more complicated. Firstly, we assume that $v_0\geq u$, and applying \eqref{lfun_derRker_estim} and the substitution $v-u=\eta\sqrt{1-r}$ we compute
\begin{align*}
\int_{v_0}^{\infty}\left(\partial_u\Rop(u,v)\right)^2 \,d\mu_\al(v)&\lesssim (1-r)^{-3}u^{-2\al-1}\int_{v_0}^{\infty}v^2\exp\left(-\frac{1+r}{1-r}(v-u)^2\right)\,dv\\
&\lesssim (1-r)^{-2\al-4}\int_{v_0}^{\infty}v^{2\al+3}\exp\left(-\frac{(v-u)^2}{1-r}\right)\,dv\\
&\lesssim (1-r)^{-2\al-7/2}\int_{0}^{\infty}((\sqrt{1-r}\eta)^{2\al+3}+u^{2\al+3})e^{-\eta^2}\,d\eta\\
&\lesssim (1-r)^{-\al-2}.
\end{align*}

Now we assume that $v_0\leq u$. This time we use \eqref{lfun_derKer_cancellation_estim}. We need to estimate two underlying components. Firstly, we have
\begin{align*}
\int_{v_0}^{\infty}v^2|R_r^{\al+1}(u,v)|^2 \,d\mu_\al(v)&\lesssim (1-r)^{-1}u^{-2\al-3}\int_{v_0}^{\infty}\exp\left(-\frac{1+r}{1-r}(v-u)^2\right)\,dv\\
&\lesssim (1-r)^{-\al-5/2}\int_{-\infty}^{\infty}\exp\left(-\frac{1+r}{1-r}v^2\right)\,dv\\
&\lesssim (1-r)^{-\al-2}.
\end{align*}
Secondly, subtly estimating and substituting $\eta=v/\sqrt{1-r}$ we obtain
\begin{align*}
&\int_{v_0}^{\infty}(1-r)^{-2}(2\sqrt{r}v-(1+r)u)^2|\Rop(u,v)|^2 \,d\mu_\al(v)\\
&\lesssim (1-r)^{-3}u^{-2\al-1}\int_{v_0}^{\infty}\big( (v-u)^2 +u^2(1-r)^4\big)\exp\left(-\frac{1+r}{1-r}(v-u)^2-\frac{2uv(1-r)}{(1+\sqrt{r})^2}\right)\,dv\\
&\lesssim (1-r)^{-\al-7/2}\int_{v_0-u}^{\infty}\big( v^2 +u^2(1-r)^4\big)\exp\left(-\frac{v^2}{1-r}-\frac{(1-r)u^2}{(1+\sqrt{r})^2}\right)\,dv\\
&\lesssim (1-r)^{-\al-3}\int_{-\infty}^{\infty}\big( (1-r)\eta^2 +(1-r)^3\big)e^{-\eta^2}\,d\eta\\
&\lesssim (1-r)^{-\al-2}.
\end{align*}
This finishes the proof of the proposition.

\end{proof}

{\thm\label{lfun_H1} For $\al\in [-1/2,\infty)^d$ there is
\begin{equation*}
\sum\limits_{n\in\NN^d}\frac{|\langle f,\lfun\rangle|}{(\ven+1)^{d+|\al|/2}}\lesssim \Vert f\Vert_{H^1(\RR_+^d,\,\mu_\al)},
\end{equation*}
uniformly in $f\in H^1(\RR_+^d,\mu_\al)$. The result is sharp in the sense that for any $\ve>0$ there exists $f\in H^1(\RR_+^d,\mu_\al)$ such that
\begin{equation*}
\sum\limits_{n\in\NN^d}\frac{|\langle f,\lfun\rangle|}{(\ven+1)^{d+|\al|/2-\ve}}=\infty.
\end{equation*}
}
\begin{proof}
We easily see that since $\Rop(x,y)$ are the tensor products of the one-dimensional kernels, we have by Proposition \ref{lfun_derRker_prop} and Lemma \ref{lfun_Rker_lemma}
\begin{align*}
\sup_{x\in\RR_+^d} \big\Vert |\nabla_x\Rop(x,\cdot)|\big\Vert_{L^2(\RR_+^d,\,\mu_\al)}&\leq\sum_{j=1}^d\Vert \partial_{x_j}R^{\alpha_j}_r(x_j,\boldsymbol{\cdot})\Vert_{L^2(\RR_+,\,\mu_{\al_j})}\prod_{i\neq j} \Vert \Ropi(x_i,\boldsymbol{\cdot})\Vert_{L^2(\RR_+,\,\mu_{\al_i})}\\
&\lesssim (1-r)^{-(|\al|+d+1)/2}, 
\end{align*}
uniformly in $r\in(0,1)$. Hence, the exponent $\gamma$ appearing in \eqref{cond C} is equal to $(|\al|+d+1)/2$. Moreover, (see for example \cite[Appendix~1]{AnkerDziubanski}) for the measure $\mu_\al$ the parameter $N$ from \ref{A2} is equal to $2|\al|+2d$. Thus, assumptions \ref{A1}-\ref{A3} are satisfied in the considered setting. Hence, by Theorem \ref{H1_general} we obtain the first part of the claim.

In order to prove sharpness, given $\ve>0$ it suffices, for any fixed $K\in\NN$, to construct an $H^1(\RR^d_+,d\mu_\al)$-atom $a_{\ell^\al}$ such that
\begin{equation}\label{lfun_counterex}
\sum\limits_{n\in\NN_+^d}\frac{|\langle a_{\ell^\al},\lfun\rangle|}{\ven^{d+|\al|/2-\ve}}\gtrsim K^\ve,
\end{equation}
where the underlying constant does not depend on $K$.

Firstly, we consider the one-dimensional case. For fixed $K$ and $2^{-\frac{1}{2\al+2}}<\delta<1$ we define
$$a_{\ell^\al}(u)=\left\{\begin{array}{ll}
(\delta^{-2\al-2}-1) (\sqrt{K}/c)^{2\al+2},& u\in(0,c \delta K^{-1/2}]=:B_1,\\
-(\sqrt{K}/c)^{2\al+2},& u\in(c\delta K^{-1/2} ,c K^{-1/2})=:B_2.
\end{array}\right. $$
The constant $c>0$ depends only on $\al$ and emerges from the estimate (compare \cite[pp.~435,~453)]{Muckenhoupt})
\begin{equation}\label{lfun_estim_small}
\lfunk(u)\simeq k^{\al/2},
\end{equation}
uniform in $k\in\NN_+$ and $0<u<ck^{-1/2}$.

A straightforward computation shows that $\int_B a_{\ell^\al}(u)\,d\mu_\al(u)=0$, where $B=B_1\cup B_2=(0,cK^{-1/2})$. Moreover, note that $\mu_\al(B)=\frac{1}{2\al+2}c^{2\al+2}K^{-\al-1}$. It can also be checked that
$$\Vert a_{\ell^\al}\Vert_{L^2(\RR^d_+,\,\mu_\al)}=(2\al+2)^{-1/2} (\sqrt{K}/c)^{\al+1}(\delta^{-2\al-2}-1)^{1/2}\leq \mu_\al(B)^{-1/2}. $$
Thus, $a_{\ell^\al}$ is an $H^1(\RR_+,\mu_\al)$-atom.

By \eqref{lfun_derivative_formula} and \eqref{lfun_estim_small}, for $0<u<ck^{-1/2}$, we have
$$-\frac{d\lfunk}{du}(u)\gtrsim u k^{\al/2}+uk^{1+\al/2}\gtrsim uk^{1+\al/2}. $$
Hence, using the mean value theorem we obtain for $k\leq K$
\begin{align*}
&\int_B a_{\ell^\al}(u)\lfunk(u)\,d\mu_\al (u)=\int_B a_{\ell^\al}(u) (u-\delta cK^{-1/2}) \frac{d\lfunk}{du}(\xi_u)\,d\mu_\al(u)\\
&=(\sqrt{K}/c)^{2\al+2}\int_B \Big((\delta^{-2\al-2}-1)\ind_{B_1}(u)+\ind_{B_2}(u) \Big)\big|u-\delta cK^{-1/2}\big|\Big(-\frac{d\lfunk}{du}(\xi_u)\Big)\,d\mu_\al (u)\\
&\gtrsim (\sqrt{K}/c)^{2\al+2} k^{1+\al/2}(\delta^{-2\al-2}-1) \int_0^{\delta c K^{-1/2}}(\delta cK^{-1/2}-u) u^{2\al+2}\,du\\
&=c^2 K^{-1}k^{1+\al/2}(\delta^{-2\al-2}-1)\delta^{2\al+4}\big((2\al+3)(2\al+4)\big)^{-1}\\
&\gtrsim K^{-1}k^{1+\al/2},
\end{align*}
where $\xi_u$ is between $u$ and $\delta c K^{-1/2}$. Finally,
$$\sum\limits_{k=1}^K \frac{|\langle a_{\ell^\al},\lfunk\rangle|}{k^{1+\al/2-\ve}}\gtrsim K^{-1}\sum_{k=1}^K k^\ve \simeq K^\ve, $$
which finishes the justification of the one-dimensional version of \eqref{lfun_counterex}.

In the multi-dimensional case we consider
$$\bold{a}_{\ell^\al}(x)=\prod_{i=1}^d a_{\ell^{\al_i}}(x_i), $$
where $a_{\ell^{\al_i}}$ is as above. It is obvious that $\mathrm{supp}\ \bold{a}_{\ell^\al}\subset B:=B(\bold{c/2\sqrt{K}},\sqrt{d}c/2\sqrt{K})$, where $\bold{c/2\sqrt{K}}=(c/2\sqrt{K},\ldots,c/2\sqrt{K})$. Moreover, $\int_B \bold{a}_{\ell^\al}(x)\,d\mu_{\al}(x)=0$. Lastly,
\begin{align*}
\Vert\bold{a}_{\ell^\al}\Vert_{L^2(\RR^d_+,\,d\mu_\al)}&=(\sqrt{K}/c)^{|\al|+d} \prod_{i=1}^d(2\al_i+2)^{-1/2}(\delta^{-2\al_i-2}-1)^{1/2}\\
&\leq (\sqrt{K}/c)^{|\al|+d} d^{-(|\al|+d)/2} \prod_{i=1}^d(2\al_i+2)^{1/2}\\
&\leq \mu_\al(B)^{-1/2},
\end{align*}
for $\delta$ sufficiently close to $1$, namely $\delta\geq \max\{ (1+d^{-\al_i-1})^{-\frac{1}{2\al_i+2}}\colon i=1,\ldots,d\}$. Hence, $\bold{a}_{\ell^\al}$ is an $H^1(\RR^d_+,\mu_\al)$-atom. 

To justify \eqref{lfun_counterex} we compute
$$\sum_{n\in\NN^d_+}\frac{|\langle \bold{a}_{\ell^\al},\lfun\rangle|}{\ven^{d+|\al|/2-\ve}}\geq K^{-d}\sum_{n\in\{1,\ldots, K\}^d}\frac{\prod_{i=1}^d n_i^{1+\al_i/2}}{\ven^{d+|\al|/2-\ve}}\gtrsim K^{-2d-|\al|/2+\ve}\prod_{i=1}^d \sum_{n_i=1}^K n_i^{1+\al_i/2}\simeq K^{\ve}.  $$
This finishes the proof of the theorem.
\end{proof}

\subsection{$L^1$ result for Laguerre functions of convolution type}

{\thm For $\al\in[-1/2,\infty)^d$ and any $\ve>0$ there is
\begin{equation*}
\sum\limits_{n\in\NN^d}\frac{|\langle f,\lfun\rangle|}{(\ven+1)^{d+|\al|/2+\ve}}\lesssim \Vert f\Vert_{L^1(\RR_+^d,\,d\mu_\al)},
\end{equation*}
uniformly in $f\in L^1(\RR_+^d,d\mu_\al)$. The result is sharp in the sense that there exists $f\in L^1(\RR_+^d,d\mu_\al)$ such that
\begin{equation}\label{lfun_L1_false}
\sum\limits_{n\in\NN^d}\frac{|\langle f,\lfun\rangle|}{(\ven+1)^{d+|\al|/2}}=\infty.
\end{equation}
}
\begin{proof}
The first part follows from \eqref{lfun_asympt} and a computation similar to the one conducted in the proof of the first part of \cite[Theorem~5.1]{Plewa}.

For the second part we assume a contrario that the sum in \eqref{lfun_L1_false} is finite for every $f\in L^1(\RR^d_+,d\mu_\al)$. Then the uniform boundedness principle and \cite[Lemma~1]{Kanjin2} yield
$$\sum_{n\in\NN^d_+} \frac{|\lfun(x)|}{\ven^{d+|\al|/2}}\lesssim 1,\qquad x\in\RR_+^d.$$ 
But this is not true. Indeed, for fixed large $K\in\NN$ and $x\in (0,c(dK)^{-1/2})^d$ we shall show that
\begin{equation}\label{lfun_L1_counter_claim}
\sum_{n\in\NN_+^d} \frac{|\lfun(x)|}{\ven^{d+|\al|/2}}\gtrsim \sum_{n\in\{1,\ldots,K\}^d} \frac{n^{\al/2}}{\ven^{d+|\al|/2}} \gtrsim \log K.
\end{equation}
For this purpose we apply the induction. If $d=1$, then by \eqref{lfun_estim_small} we have
$$\sum_{k=1}^K \frac{|\lfunk(u)|}{k^{1+\al/2}}\gtrsim \sum_{k=1}^K \frac{k^{\al/2}}{k^{1+\al/2}}\gtrsim\log K,\qquad 0<u<cK^{-1/2}. $$
Now we assume that \eqref{lfun_L1_counter_claim} holds for some $d\geq 1$ and will justify it for $d+1$. We estimate for $\beta\geq -1/2$
\begin{align*}
\sum_{n\in\{1,\ldots,K\}^d}\sum_{k=1}^K \frac{|\ell_{(n,k)}^{(\al,\beta)}(x,u)|}{(\ven+k)^{d+1+|\al|/2+\beta/2}}&\gtrsim \sum_{n\in\{1,\ldots,K\}^d}\!\!\!\!\! n^{\al/2}\!\!\!\!\!\sum_{\ven/d\leq k\leq K} \Big(\frac{k}{\ven+k}\Big)^{\beta/2}(\ven+k)^{-d-1-|\al|/2}\\
&\gtrsim \sum_{n\in\{1,\ldots,K\}^d}n^{\al/2}\big((2\ven)^{-d-|\al|/2}-(\ven+K)^{-d-|\al|/2}\big)\\
&\gtrsim \log K - 1,
\end{align*}
and this finishes the proof of the theorem.

\end{proof}

\section{Standard Laguerre setting}\label{S4}

The {\it standard Laguerre functions} $\{\Lfunk\}_{k\in\mathbb{N}}$ of order $\al>-1$ are defined on $\RR_+$ by
$$\Lfunk(u)=2^{-1/2} u^{\al/2}\lfunk(u^{1/2}),\qquad u>0,$$
and in the multi-dimensional case as the tensor product of the one-dimensional functions. The system $\{\Lfun\}_{n\in\NN^d}$ form an orthonormal basis in $L^2(\RR_+^d,\,dx)$.
We shall use the pointwise asymptotic estimates similar to those in \eqref{lfun_asympt} (see \cite[p.~435]{Muckenhoupt} and \cite[p.~699]{AskeyWainger})
\begin{equation}\label{Lfun_asympt}
\vert \Lfunk(u)\vert\lesssim\left\{ \begin{array}{ll}
(u\nu)^{\al/2},  & 0<u\leq 1/\nu,\\
(u\nu)^{-1/4}, & 1/\nu<u\leq\nu/2,\\
(\nu(\nu^{1/3}+\vert u- \nu\vert))^{-1/4}, & \nu/2<u\leq 3\nu/2,\\
\exp(-\gamma u), & 3\nu/2<u<\infty,
\end{array}\right.
\end{equation}
where $\gamma>0$ depends only on $\alpha$. 

The above estimates and formula
$$(\Lfunk)'(u)=-k^{1/2}\mathcal{L}_{k-1}^{\al+1}(u)u^{-1/2}-\frac{1}{2}\Lfunk(u)+\frac{\alpha}{2u}\Lfunk(u), $$
where $\mathcal{L}_{-1}^{\al+1}\equiv 0$, yield 
\begin{equation}\label{Lfun_bound}
\Vert \Lfunk\Vert_{L^\infty(\RR_+)} \lesssim 1, \quad\al\geq 0,\qquad \Vert (\Lfunk)'\Vert_{L^\infty(\RR_+)} \lesssim k+1, \quad \al\geq 2.\\
\end{equation}
For the first estimate compare \cite[p.~94]{StempakTohoku}.

\subsection{$L^1$ result for the standard Laguerre functions}

Before we shall investigate Hardy's inequality associated with the standard Laguerre functions we will study its $L^1$-analogue. The following lemma will be needed. 

{\lm\label{Lfun_alfa_lem} If $\al\in[0,\infty)^d\setminus\{\bold{0}\}$, then
$$\sum\limits_{n\in\NN^d}\frac{|\langle f,\Lfun\rangle|}{(\ven+1)^d}\lesssim \Vert f\Vert_{L^1(\RR_+^d)}+\sum\limits_{n\in\NN^{d}}\frac{|\langle f,\mathcal{L}_n^{\al+\bold{2}}\rangle|}{\prod_{i=1}^d(n_i+1)},$$
uniformly in $f\in L^1(\RR_+^d)$.
}
\begin{proof}
We will use the induction on $d$. For $d=1$ the result is proved in \cite[p.~335]{Kanjin1}, where the main tool is the estimate 
\begin{equation*}
|\langle g,\mathcal{L}_k^\beta\rangle|\lesssim |\langle g,\mathcal{L}_{k-1}^{\beta+2}\rangle|+\sum_{j=k}^\infty |\langle g,\mathcal{L}_j^{\beta+2}\rangle|(k/j)^{\beta/2}j^{-1},\qquad g\in L^1(\RR_+),\ k\in\NN_+,
\end{equation*}
see \cite[p.~401]{Askey}. Here $\beta\geq 0$ (although in the referred estimate the parameter is assumed to be non-zero, the inequality holds for $\beta=0$ as well). 

We assume that the claim holds in a dimension $d\geq1$, and we will prove it in the dimension $d+1$. Fix $f\in L^1(\RR_+^{d+1})$ and $(\al,\beta)\in[0,\infty)^{d+1}\setminus\{\bold{0}\}$. Without any loss of generality we may assume that $\al\in[0,\infty)^d\setminus\{\bold{0}\}$ and $\beta\geq 0$. For $n\in\NN^d$ we denote $f_n^\al(u)=\langle f(\cdot,u),\Lfun\rangle$, where $(x,u)$ are the coordinates in $\RR^{d+1}_+$. Clearly, $f_n^\al\in L^1(\RR_+)$. Using the tool invoked above we obtain for $\beta\geq 0$
\begin{equation*}
|\langle f_n^\al,\mathcal{L}_k^\beta\rangle|\lesssim |\langle f_n^\al,\mathcal{L}_{k-1}^{\beta+2}\rangle|+\sum_{j=k}^\infty |\langle f_n^\al,\mathcal{L}_j^{\beta+2}\rangle|(k/j)^{\beta/2}j^{-1},\qquad k\in\NN_+,\ n\in\NN^d.
\end{equation*} 
Hence,
\begin{align*}
&\sum\limits_{n\in\NN^d}\sum_{k=0}^\infty\frac{|\langle f,\mathcal{L}_{(n,k)}^{(\al,\beta)}\rangle|}{(\ven+k+1)^{d+1}}\\
&\lesssim \sum\limits_{n\in\NN^d}\bigg(\frac{|\langle f_n^\al,\mathcal{L}_0^\beta\rangle|}{(\ven+1)^{d+1}}+ \sum_{k=1}^\infty \frac{|\langle f_n^\al,\mathcal{L}_{k-1}^{\beta+2}\rangle|+\sum_{j=k}^\infty |\langle f_n^\al,\mathcal{L}_{j}^{\beta+2}\rangle| (k/j)^{\beta/2}j^{-1}}{(\ven+k+1)^{d+1}}\bigg)\\
&\lesssim \sum\limits_{n\in\NN^d}\bigg(\frac{\Vert f_n^\al\Vert_{L^1(\RR+)}}{(\ven+1)^{d+1}}+ \sum_{k=1}^\infty\frac{|\langle f_n^\al,\mathcal{L}_{k-1}^{\beta+2}\rangle|+\sum_{j=k}^\infty |\langle f_n^\al,\mathcal{L}_{j}^{\beta+2}\rangle|j^{-1}}{(\ven+k+1)^{d+1}}\bigg)\\
&\lesssim \Vert f\Vert_{L^1(\RR_+^d)}+\sum\limits_{n\in\NN^d,\, k\in\NN}\frac{|\langle f,\mathcal{L}_{(n,k)}^{(\al,\beta+2)}\rangle|}{(\ven+k+1)^{d+1}}+\sum_{n\in\NN^d}\sum\limits_{j=1}^\infty |\langle f_n^\al,\mathcal{L}_{j}^{\beta+2}\rangle|j^{-1}\sum\limits_{k=1}^j (\ven+k+1)^{-d-1}.
\end{align*}

Firstly, we denote $f_k^{\beta+2}(x)=\langle f(x,\cdot),\mathcal{L}_k^{\beta+2}\rangle$. Clearly, $f_k^{\beta+2}\in L^1(\RR^d_+)$. Then we estimate the first sum from the last line above using the inductive hypothesis
\begin{align}
\sum\limits_{n\in\NN^d,\, k\in\NN}\frac{|\langle f,\mathcal{L}_{(n,k)}^{(\al,\beta+2)}\rangle|}{(\ven+k+1)^{d+1}}&\leq\sum_{k=0}^\infty \frac{1}{k+1}\sum_{n\in\NN^d}\frac{|\langle f_k^{\beta+2},\Lfun\rangle|}{(\ven+1)^d}\label{Lfun_alfa_lemma_1}\\
&\lesssim \sum_{k=0}^\infty \frac{1}{k+1}\Big(\Vert f_k^{\beta+2}\Vert_{L^1(\RR_+^d)}+\sum_{n\in\NN^d}\frac{|\langle f_k^{\beta+2},\mathcal{L}_n^{\al+\bold{2}}\rangle|}{\prod_{i=1}^d(n_i+1)}\Big)\nonumber\\
&\lesssim \int_{\RR_+^d} \bigg(\sum_{k=0}^\infty\frac{|\langle f(x,\cdot),\mathcal{L}_k^{\beta+2}\rangle|}{k+1}\bigg)\,dx+\sum_{n\in\NN^d,\,k\in\NN}\frac{|\langle f,\mathcal{L}_{(n,k)}^{(\al+\bold{2},\beta+2)}\rangle|}{(k+1)\prod_{i=1}^d(n_i+1)}.\nonumber
\end{align}
Note that the first summand is bounded by $\Vert f\Vert_{L^1(\RR_+^{d+1})}$, since for $\delta>0$
\begin{equation}\label{Lfun_sum_estim_d=1}
\sum_{k=0}^\infty \frac{|\mathcal{L}_k^{\delta}(u)|}{k+1}\lesssim 1,\qquad u>0,
\end{equation}
see \cite[Lemma~3]{Kanjin2}.

In order to estimate the second remaining sum we use the bound
$$\sum_{k=1}^\infty (\ven+k+1)^{-d-1}\lesssim (\ven+1)^{-d},\qquad n\in\NN^d,$$
to obtain
$$\sum_{n\in\NN^d}\sum\limits_{j=1}^\infty |\langle f_n^\al,\mathcal{L}_{j}^{\beta+2}\rangle|j^{-1}\sum\limits_{k=1}^j (\ven+k+1)^{-d-1}\lesssim \sum_{j=0}^\infty \frac{1}{j+1} \sum_{n\in\NN^d}\frac{|\langle f_j^{\beta+2},\Lfun\rangle|}{(\ven+1)^d},$$
and the same quantity was estimated in \eqref{Lfun_alfa_lemma_1}. This concludes the proof of the lemma.
\end{proof}

{\thm\label{Lfun_L1} For $\al\in[0,\infty)^d\setminus\{\bold{0}\}$ there is
\begin{equation*}
\sum\limits_{n\in\NN^d}\frac{|\langle f,\Lfun\rangle|}{(\ven+1)^{d}}\lesssim \Vert f\Vert_{L^1(\RR_+^d)},
\end{equation*}
uniformly in $f\in L^1(\RR_+^d)$. The result is sharp in the sense that for any $\ve>0$ there exists $f\in L^1(\RR_+^d)$ such that
\begin{equation*}
\sum\limits_{n\in\NN^d}\frac{|\langle f,\Lfun\rangle|}{(\ven+1)^{d-\ve}}=\infty.
\end{equation*}
}
We remark that in the admissible range of $\al$'s we cannot include $\bold{0}$. In that case the corresponding series would also diverge with the exponent equal to $d$ (compare \cite[Proposition]{Kanjin2}).
\begin{proof}
Firstly, note that sharpness follows from Theorem \ref{Lfun_H1}, which is stated and proved below.

In order to prove the remaining claim we apply Lemma \ref{Lfun_alfa_lem}. It suffices to justify that
$$\sum\limits_{n\in\NN^d}\frac{|\mathcal{L}_n^{\al+\bold{2}}(x)|}{\prod_{i=1}^d(n_i+1)}\lesssim 1,  $$
uniformly in $x\in\RR_+^d$. Note that
$$\sum\limits_{n\in\NN^d}\frac{|\mathcal{L}_n^{\al+\bold{2}}(x)|}{\prod_{i=1}^d(n_i+1)}= \prod_{i=1}^d\sum\limits_{n_i=0}^\infty \frac{|\mathcal{L}_{n_i}^{\al_i+2}(x_i)|}{(n_i+1)}. $$
Hence, \eqref{Lfun_sum_estim_d=1} implies the claim.

%Now let $\al\in[0,\infty)^d\setminus\{\bold{0}\}$. In this case we apply Lemma \ref{Lfun_alfa_lem}. It suffices to justify that
%$$\sum\limits_{n\in\NN^{d}}\frac{|\langle f,\mathcal{L}_n^{\al+\bold{2}}\rangle|}{\prod_{i=1}^d(n_i+1)}\lesssim \Vert f\Vert_{L^1(\RR_+^d)}.$$
%This is true since \eqref{Lfun_sum_estim_d=1} yields
%$$\sum_{n\in\NN^d}\frac{|\mathcal{L}_n^{\al+\bold{2}}(x)|}{\prod_{i=1}^d(n_i+1)} $$
%
%$$\langle f,\mathcal{L}_n^{\al+\bold{2}}\rangle=\langle f_{\bar{n}}^{\bar{\al}+\bold{2}},\mathcal{L}_{n_d}^{\al_d+2}\rangle, $$
%where $n=(\bar{n},n_d)$, $\al=(\bar{\al},\al_d)$, and $f_{\bar{n}}^{\bar{\al}+\bold{2}}(u)=\langle f(\cdot,u),\mathcal{L}_{\bar{n}}^{\bar{\al}+\bold{2}}\rangle\in L^1(\RR_+)$. Moreover, \eqref{Lfun_bound} implies that $\Vert f_{\bar{n}}^{\bar{\al}+\bold{2}}\Vert_{L^1(\RR_+)}\leq \Vert f\Vert_{L^1(\RR^d_+)}.$ Hence, it suffices to iterate the already proved one-dimensional version of the theorem.
\end{proof}

\subsection{$H^1$ result for standard Laguerre functions}

Theorem \ref{Lfun_L1} implies that Hardy's inequality for standard Laguerre functions holds for $\al\in[0,\infty)^d\setminus\{\bold{0}\}$ with the admissible exponent $E=d$. In fact, in Theorem \ref{Lfun_H1} we show that the exponent is sharp. Moreover, Hardy's inequality is also valid for $\al=\bold{0}$. For this purpose we shall estimate the kernels associated with the family of operators $\{\Rop\}_{r\in(0,1)}$ corresponding to the functions $\{\Lfun\}_{n\in\NN^d}$, similarly as we did in Section \ref{S3}. We could restrict the reasoning to $\al=\bold{0}$, but at a low cost we obtain the auxiliary results for more general range of $\al$'s.

On the other hand, we could use \ref{A3} instead of easier \eqref{cond C} in order to fill the gap in the admissible range of the Laguerre type parameter in Proposition \ref{Lfun_derRker_prop}. However, Lemma \ref{Lfun_alfa_lem} covers this gap in Theorem \ref{Lfun_H1}, and therefore we do not care about the restraint in Proposition \ref{Lfun_derRker_prop}.

The family of operators $\{\Rop\}_{r\in(0,1)}$ associated with the functions $\Lfun$ is also a family of integral operators. Similarly to Section \ref{S3}, the kernels associated with the operators $\Rop$ are defined by
$$\Rop(x,y)=\sum_{n\in\NN^d}^\infty r^{\ven}\Lfun(x)\Lfun(y),\qquad x,y\in\RR^d_+,\ r\in(0,1). $$
The one-dimensional kernels have the explicit representation (compare \cite[p.~102]{Szego})
\begin{equation}\label{Lfun_Rker_explic}
\Rop(u,v)=(1-r)^{-1}r^{-\al/2}\exp\Big(-\frac{1}{2}\frac{1+r}{1-r}(u+v)\Big)I_{\al }\Big(\frac{2r^{1/2}}{1-r}\sqrt{uv}\Big).
\end{equation}

We shall estimate the kernels $\Rop(u,v)$. Using \eqref{Bessel_estim} we obtain
\begin{align}
&|\Rop(u,v)|\nonumber\\
&\lesssim\left\{\begin{array}{ll}
(1-r)^{-\al-1}(uv)^{\al/2}\exp\big(-\frac{1}{2}\frac{1+r}{1-r}(u+v) \big), & v\leq \frac{(1-r)^2}{4ru},\\
(1-r)^{-1/2}r^{-\al/2-1/4}(uv)^{-1/4}\exp\big(-\frac{1}{2}\frac{1+r}{1-r}(\sqrt{v}-\sqrt{u})^2-\frac{\sqrt{uv}(1-r)}{(1+\sqrt{r})^2} \big), & v\geq \frac{(1-r)^2}{4ru}.\label{Lfun_Rker_estim}
\end{array}\right.
\end{align}

{\lm\label{Lfun_Rker_standard} For $\al\geq 0$ there is
$$\sup_{u\in\RR_+}\Vert \Rop(u,\boldsymbol{\cdot})\Vert_{L^2(\RR_+)}\lesssim (1-r)^{-1/2},\qquad r\in(0,1). $$}
\begin{proof}
For $0<r\leq 1/2$, Parseval's identity and \eqref{Lfun_bound} imply
$$\sup_{u\in\RR_+}\Vert \Rop(u,\boldsymbol{\cdot})\Vert_{L^2(\RR_+)}\leq  \Big(\sum_{k=0}^\infty 2^{-2 k}\Vert \Lfunk\Vert_{L^{\infty}(\RR_+)}^2\Big)^{1/2}\lesssim 1.$$

For $1/2<r<1$, we denote $v_0=(1-r)^2/(4ru)$, and estimate the integrals over $(0,v_0]$ and $(v_0,\infty)$. Using the substitution $\eta=v(1+r)/(1-r)$ we obtain
\begin{align*}
\int_0^{v_0} \Rop(u,v)^2\,dv&\lesssim  (1-r)^{-2\al-2}\int_0^{v_0}(uv)^\al\exp\Big(-\frac{1+r}{1-r}v\Big)\,dv\\
&\lesssim (1-r)^{-2\al-2}\int_0^{\infty}(1-r)^{2\al}e^{-\eta}(1-r)\,d\eta\\
&\lesssim (1-r)^{-1}.
\end{align*}
Similarly, we substitute $\eta=\sqrt{v}$, and for $u\geq v_0$ we get
\begin{align*}
\int_{v_0}^{\infty} |\Rop(u,v)|^2\,dy&\lesssim  (1-r)^{-1}\int_{v_0}^{\infty}(uv)^{-1/2}\exp\Big(-\frac{1+r}{1-r}(\sqrt v-\sqrt u)^2\Big)\,dv\\
&\lesssim  (1-r)^{-1}\int_{-\infty}^{\infty}(1-r)^{-1/2}\exp\Big(-\frac{1+r}{1-r}(\eta-\sqrt{u})^2\Big)\,d\eta\\
&\lesssim (1-r)^{-1},
\end{align*}
and for $u\leq v_0$ we use $\eta=\sqrt{v}-\sqrt{u}$ and receive 
\begin{align*}
\int_{v_0}^{\infty} |\Rop(u,v)|^2\,dy&\lesssim  (1-r)^{-1}\int_{0}^{\infty}(1-r)^{-1}(\eta+\sqrt{u})\exp\Big(-\frac{1+r}{1-r}\eta^2\Big)\,d\eta\\
&\lesssim (1-r)^{-1}.
\end{align*}
Combining the above gives the claim.
\end{proof}

Using the same formula for the derivatives of the Bessel function $I_\al$ as in \eqref{lfun_derKer_long}, we obtain

\begin{equation}\label{Lfun_derKer_formula}
\partial_u\Rop(u,v)=\frac{r\sqrt{v}}{\sqrt{u}(1-r)}R^{\alpha+1}_r(u,v)+\left(\frac{\al}{2u}-\frac{1+r}{2(1-r)}\right)\Rop(u,v).
\end{equation}

Hence, the above equality and \eqref{Lfun_Rker_estim} yield
\begin{align*}
&\Big|\partial_u\Rop(u,v)\Big|^2\\
&\lesssim\left\{\begin{array}{ll}
(1-r)^{-2-2\al}(uv)^\al\exp\left(-\frac{1+r}{1-r}(u+v) \right)(u^{-2}+(1-r)^{-2}), & v\leq \frac{(1-r)^2}{4ru},\\
(1-r)^{-3}r^{-\al-1} (uv)^{-1/2}\exp\left(-\frac{1+r}{1-r}(\sqrt{v}-\sqrt{u})^2-\frac{2\sqrt{uv}(1-r)}{(1+\sqrt{r})^2} \right)(1+v/u), & v\geq \frac{(1-r)^2}{4ru}.
\end{array}\right.
\end{align*}

{\prop\label{Lfun_derRker_prop} For $\al\geq 2$ or $\al=0$ we have
$$\sup_{u\in\RR_+}\left\Vert\partial_u\Rop(u,\boldsymbol{\cdot})\right\Vert_{L^2(\RR_+)}\lesssim (1-r)^{-3/2},\qquad r\in(0,1). $$}
\begin{proof}
If $0<r\leq 1/2$, then we again apply Parseval's identity and \eqref{Lfun_bound} obtaining
\begin{equation*}
\left\Vert\partial_u\Rop(u,\boldsymbol{\cdot})\right\Vert_{L^2(\RR_+)}\leq\Big(\sum_{k=0}^\infty 2^{-2k}\Vert(\Lfunk)'\Vert_{L^{\infty}(\RR_+)}^2\Big)^{1/2}\lesssim 1.
\end{equation*}

From now on we assume that $1/2<r<1$. Firstly we deal with the case $\al\geq 2$. We use the notation $v_0=(1-r)^2/(4ru)$ again and split the integration over two intervals: $(0,v_0]$ and $(v_0,\infty)$. In the first case, we estimate
\begin{align*}
&\int_0^{v_0}|\partial_u\Rop(u,v)|^2\,dv\\
 &\lesssim (1-r)^{-2\al-2}\int_0^{y_0}(u^{\al-2}+u^{\al}(1-r)^{-2})\exp{\left(-\frac{1+r}{1-r}(v+u)\right)}v^\al\, dv\\
&\lesssim (1-r)^{-3}.
\end{align*} 

In the second case we firstly assume that $v_0\leq u$. Applying the pointwise estimate of $\big|\partial_u\Rop(u,v)\big|$, and substituting $\eta=\sqrt{v}$, we compute
\begin{align*}
\int_{v_0}^{\infty}|\partial_u\Rop(u,v)|^2 \,dv&\lesssim (1-r)^{-3}u^{-1/2}\int_{0}^{\infty}\exp\left(-\frac{1+r}{1-r}(\eta-\sqrt{u})^2\right)(1+\eta^2u^{-1})\,d\eta\\
&\lesssim (1-r)^{-7/2}\int_{-\infty}^{\infty}\exp\left(-\frac{1+r}{1-r}\eta^2\right)(1+\eta^2(1-r)^{-1})\,d\eta\\
&\lesssim (1-r)^{-3}.
\end{align*}

For $u\leq v_0$ we use $\eta=\sqrt{v}-\sqrt{u}$ and receive 
\begin{align*}
\int_{v_0}^{\infty} |\partial_u\Rop(u,v)|^2\,dv&\lesssim (1-r)^{-4}\int_{v_0}^{\infty} (1+v^2(1-r)^{-2})\exp\Big(-\frac{1+r}{1-r}(\sqrt{v}-\sqrt{u})^2\Big)\,dv\\
&\lesssim (1-r)^{-6}\int_{v_0}^{\infty}((1-r)^2+v^2) \exp\Big(-\frac{1+r}{1-r}(\sqrt{v}-\sqrt{u})^2\Big)\,dv\\
&\lesssim  (1-r)^{-6}\int_{0}^{\infty}((1-r)^2+(\eta+\sqrt{u})^4)(\eta+\sqrt{u})\exp\Big(-\frac{1+r}{1-r}\eta^2\Big)\,d\eta\\
&\lesssim (1-r)^{-3}.
\end{align*}
This finishes the proof for $\al\geq 2$. Note that the estimates for $v\in[v_0,\infty)$ are also valid for $\al=0$. Therefore, it suffices to take care of the case of integration over $(0,v_0)$. Note that \eqref{Lfun_Rker_estim} and \eqref{Lfun_derKer_formula} imply
$$|\partial_u R^0_r(u,v)|^2\lesssim (1-r)^{-4}\exp\Big(-\frac{1+r}{1-r}v \Big)(1+v^2(1-r)^{-2}),\qquad  v\leq \frac{(1-r)^2}{4ru}. $$
Now it is easily seen that substituting $\eta=v(1+r)(1-r)^{-1}$, we have
\begin{align*}
\int_0^{v_0}|\partial_u R^0_r(u,v)|^2\,dy\lesssim (1-r)^{-4}\int_0^{\infty}(1+\eta^2)e^{-\eta}(1-r)\, d\eta\lesssim (1-r)^{-3},
\end{align*}
and it finishes the proof of the proposition.
\end{proof}

Now we calculate the parameter $\gamma$ in \eqref{cond C} (multi-dimensional case). Combining Lemma \ref{Lfun_Rker_standard} and Proposition \ref{Lfun_derRker_prop} we receive
$$\big\Vert | \nabla_x\Rop(x,\boldsymbol{\cdot})|\big\Vert_{L^2(\RR_+^d)}\lesssim (1-r)^{-(d+2)/2},$$
uniformly in $x\in\RR_+^d$ and $r\in(0,1)$. Hence, for the standard Laguerre functions we obtained $\gamma=1+d/2$.

{\thm\label{Lfun_H1} For $\al\in[0,\infty)^d$ there holds 
\begin{equation}
\sum\limits_{n\in\NN^d}\frac{|\langle f,\Lfun\rangle|}{(\ven+1)^d}\lesssim \Vert f\Vert_{H^1(\RR_+^d)},
\end{equation}
uniformly in $f\in H^1(\RR_+^d)$. The result is sharp in the sense that for any $\ve>0$ there exists $f\in H^1(\RR_+^d)$ such that
$$\sum\limits_{n\in\NN^d}\frac{|\langle f,\Lfun\rangle|}{(\ven+1)^{d-\ve}}=\infty. $$}
\begin{proof}
For $\al\in(\{0\}\cup[2,\infty))^d$ we have $\gamma=1+d/2$ in \eqref{cond C} for the standard Laguerre functions. Moreover, for Lebesgue measure the parameter $N$ from \ref{A2} is equal to $d$. Hence, \ref{A1}-\ref{A3} are satisfied, and thus we get the first part of the claim by Theorem \ref{H1_general} for the indicated $\al$'s. Then, Theorem \ref{Lfun_L1} justifies the claim for the full range of the parameter $\al$.

In order to prove the second part we consider firstly the one-dimensional case. For fixed $\ve>0$ and $K\in\NN$ we shall construct an $H^1(\RR_+)$-atom $a_{\mathcal{L}^\al}$ such that
$$\sum\limits_{k=0}^\infty\frac{|\langle a_{\mathcal{L}^\al},\Lfun\rangle|}{(k+1)^{1-\ve}}\gtrsim K^\ve,$$
with the underlying constant independent of $K$.

Firstly, let $\al>0$ and fix $\ve>0$ and $K\in\NN$. For $\delta\in (0,1/2)$ we define
$$a_{\mathcal{L}^\al}(u)=\left\{\begin{array}{ll}
\delta c^{-1}(1-\delta)^{-1} K,& u\in(c\delta K^{-1},c K^{-1}),\\
-c^{-1}K,& u\in(0,c\delta K^{-1}),
\end{array}\right. $$
where $c>0$ is a constant depending only on $\al$. It is easy to check that $a_{\mathcal{L}^\al}$ is an $H^1(\RR_+^d)$-atom.

There are known estimates for functions $\Lfunk$ (see \cite[pp.~435,~453)]{Muckenhoupt})
\begin{equation*}
A(ku)^{\al/2}\leq \Lfunk(u)\leq B(ku)^{\al/2},\qquad 0<u\leq \frac{c}{k},
\end{equation*}
where $c$ is the same constant as in the definition of $a_{\mathcal{L}^\al}$, and $A,B>0,$ depend only on $\al$. Hence, for $1\leq k\leq K$, we have
\begin{align*}
\int_{\RR_+}a_{\mathcal{L}^\al}(u)\Lfunk(u)\,du&\geq \frac{\delta AKk^{\al/2}}{c(1-\delta)}\int_{c\delta/K}^{c/K}u^{\al/2}\,du-\frac{KBk^{\al/2}}{c} \int_{0}^{c\delta/K}u^{\al/2}\,du\\
&=\frac{2k^{\al/2}A\delta c^{\al/2}}{(\al+2)(1-\delta)K^{\al/2}}\left(1-\delta^{\al/2}(\delta+B/A)+\delta^{1+\al/2}B/A  \right)\\
&\gtrsim \frac{k^{\al/2}\delta}{K^{\al/2}(1-\delta)}\left(1- \delta^{\al/2}(\delta+B/A)\right)
\end{align*}
Choosing $\delta$ sufficiently small and independently of $K$ (to be precise one can take $\delta=(A/2b)^{2/\al}$) we obtain
$$\langle a_{\mathcal{L}^\al},\Lfunk\rangle\gtrsim k^{\al/2}K^{-\al/2}.$$
Thus, we have
\begin{align*}
\sum\limits_{k=0}\infty\frac{|\langle a_{\mathcal{L}^\al},\Lfunk\rangle|}{(k+1)^{1-\ve}}\gtrsim K^{-\al/2}\sum\limits_{k=1}^K k^{\al/2+\ve-1}\simeq K^{\ve},
\end{align*}
which finishes the proof for $d=1$ and $\al>0$. If $\al=0$, then we use the atom constructed in the proof of Theorem \ref{lfun_H1}. Namely,
$$a_{\mathcal{L}^0}(u)=a_{\ell^0}(\sqrt{u}),\qquad u>0. $$
Computing directly one can show that $a_{\mathcal{L}^0}$ is an $H^1(\RR_+)$-atom. Moreover, note that for $k\leq K$ we have
$$\langle a_{\mathcal{L}^0},\mathcal{L}^0_k\rangle=\int_{\RR_+} a_{\mathcal{L}^0}(u)\mathcal{L}_k^0(u)\,du\simeq \int_{\RR_+} a_{\ell^0}(u)\ell^0(u)\,d\mu_0(u)\geq k K. $$
Hence,
\begin{align*}
\sum\limits_{k=0}^\infty\frac{|\langle a_{\mathcal{L}^0},\mathcal{L}^0_k\rangle|}{(k+1)^{1-\ve}}\gtrsim K^{-1}\sum\limits_{k=1}^K k^{\ve}\simeq K^{\ve}.
\end{align*}

 A similar construction can be performed in the multi-dimensional case for $\al\in [0,\infty)^d$, but we omit the details. 
\end{proof}

\section*{Appendix}
\begin{proof}[Proof of Theorem \ref{H1_general}]
Firstly, we shall prove the claim for atoms, i.e.
$$\sum\limits_{n\in\NN^d}\frac{|\langle a,\fun\rangle|}{(\ven+1)^E}\lesssim 1, $$
uniformly in $H^1(X,\mu)$-atoms.
We apply the idea from \cite{LiYuShi}: using the asymptotic estimate for the beta function, namely $B(k,m)\simeq \Gamma(m) k^{-m}$, for large $k$ and fixed $m$, we have
\begin{align*}
\sum_{n\in\NN^d}\frac{\vert\langle  a, \fun\rangle\vert}{(\ven+1)^{E}}&\lesssim\sum_{n\in\NN^d}\int_0^1 r^{2\ven}(1-r)^{E-1}\vert\langle a, \fun\rangle\vert \,dr\\
&\leq\int_0^1 (1-r)^{E-1}\Big(\sum_{n\in\NN^d}r^{2\ven}\Big)^{1/2}\Big(\sum_{n\in\NN^d}r^{2\ven}\vert\langle a,\fun\rangle\vert^2\Big)^{1/2}\,dr\\
&\lesssim\int_0^1 (1-r)^{E-1}(1-r)^{-d/2}\Vert R_r a\Vert_{L^2(X,\,\mu)}\,dr\\
&\lesssim 1,
\end{align*}
where in the last step we used Lemma \ref{Lemma_S2}.

In order to extend the result to whole $H^1(X,\mu)$ we need a simple continuity argument. Let us define $T(f)=\{\langle f,\fun\rangle\}_{n\in\NN^d}$ for $f\in H^1(X,\mu)$. Our aim is to prove that $T\colon H^1(X,\mu)\rightarrow \ell^1((\ven+1)^{-E})$, is bounded. Note that assumption \ref{A1} and \eqref{H1_L1_norm_ineq} yield
$$\vert\langle f,\fun\rangle\vert\leq \Vert \fun\Vert_{L^{\infty}(X,\,\mu)} \Vert f\Vert_{L^1(X,\,\mu)}\lesssim \Vert \fun\Vert_{L^\infty(X,\,\mu)}\Vert f\Vert_{H^1(X,\,\mu)},\qquad n\in\NN^d,\ f\in H^1(X,\,\mu).$$
Let us denote 
$$\omega(n)=(\ven+1)^{-d-E}(\Vert \fun\Vert_{L^\infty(X,\,\mu)}+1)^{-1}.$$
The estimate above and the inequality $E>0$ imply that  $T\colon H^1(X,\mu)\rightarrow \ell^1(\omega(n))$ is bounded. Note also that
\begin{equation}\label{ell_norms_ineq}
 \Vert \boldsymbol{\cdot}\Vert_{\ell^1(\omega(n))}\leq \Vert \boldsymbol{\cdot}\Vert_{\ell^1((\vert n\vert+1)^{-E})}.
\end{equation}

Let us fix $f\in H^1(X,\mu)$ and $f=\sum_{i=0}^{\infty}\lambda_i a_i$ be an atomic decomposition of $f$. Denote $f_m=\sum_{i=0}^{m}\lambda_i a_i$ and note that $T(f_m)$ is a Cauchy sequence in $\ell^1((\vert n\vert+1)^{-E})$. Indeed, we have for $l<m$,
$$\Vert T(f_m)-T(f_l) \Vert_{\ell^1((\vert n\vert+1)^{-E})}\leq \sum_{i=l+1}^{m}\vert\lambda_i\vert\Vert T(a_i)\Vert_{\ell^1((\ven+1)^{-E})}\lesssim \sum_{i=l+1}^m \vert \lambda_i\vert. $$
Hence, $T(f_m)$ converges to a sequence $g$ in $\ell^1((\ven+1)^{-E})$ and, by \eqref{ell_norms_ineq}, also in $\ell^1(\omega(n))$. Since $T\colon H^1(X,\mu)\rightarrow \ell^1(\omega(n))$ is bounded we have $T(f_m)\rightarrow T(f)$ in $\ell^1(\omega(n))$, therefore $g=T(f)$. To obtain the boundedness of $T\colon H^1(X,\mu)\rightarrow \ell^1((\ven+1)^{-E})$ we fix $\ve>0$, take $m$ such that $\Vert T(f-f_m)\Vert_{\ell^1((\ven+1)^{-E})}<\ve$, and calculate
\begin{align*}
\Vert T(f)\Vert_{\ell^1((\ven+1)^{-E})}&\leq \Vert T(f-f_m)\Vert_{\ell^1((\ven+1)^{-E})}+\Vert T(f_m)\Vert_{\ell^1((\ven+1)^{-E})}\\
&\leq \ve +\sum_{i=0}^m \vert \lambda_i\vert \Vert T(a_i)\Vert_{\ell^1((\ven+1)^{-E})}\\
&\lesssim \ve+\Vert f\Vert_{H^1(X,\,\mu)}.
\end{align*}
This finishes the proof.
\end{proof}

\end{document}